\documentclass[11pt]{article}

\usepackage{mathrsfs}
\usepackage{amssymb}
\usepackage{amscd}
\usepackage{amsfonts}

\font\tenmsb=msbm10
\font\sevenmsb=msbm7
\font\fivemsb=msbm5

\catcode`\@=11
\ifx\amstexloaded@\relax\catcode`\@=\active
\endinput\else\let\amstexloaded@\relax\fi
\def\spaces@{\space\space\space\space\space}
\def\spaces@@{\spaces@\spaces@\spaces@\spaces@\spaces@}
\def\space@.{\futurelet\space@\relax}
\space@.
\def\Err@#1{\errhelp\defaulthelp@\errmessage{AmS-TeX error: #1}}
\def\relaxnext@{\let\next\relax}
\def\accentfam@{7}
\def\noaccents@{\def\accentfam@{0}}
\def\Cal{\relaxnext@\ifmmode\let\next\Cal@\else
\def\next{\Err@{Use \string\Cal\space only in math mode}}\fi\next}
\def\Cal@#1{{\Cal@@{#1}}}
\def\Cal@@#1{\noaccents@\fam\tw@#1}
\def\Bbb{\relaxnext@\ifmmode\let\next\Bbb@\else
\def\next{\Err@{Use \string\Bbb\space only in math mode}}\fi\next}
\def\Bbb@#1{{\Bbb@@{#1}}}
\def\Bbb@@#1{\noaccents@\fam\msbfam#1}
\newfam\msbfam
\textfont\msbfam=\tenmsb
\scriptfont\msbfam=\sevenmsb
\scriptscriptfont\msbfam=\fivemsb

\def\NN{{\Bbb N}}
\def\ZZ{{\Bbb Z}}
\def\RR{{\Bbb R}}
\def\QQ{{\Bbb Q}}
\def\CC{{\Bbb C}}

\def\pa{{\partial}}

\def\beq{\begin{equation}}
\def\eeq{\end{equation}}

\def\skipaline{\removelastskip\vskip12pt plus 1pt minus 1pt}

\def\Proof{\removelastskip\skipaline
\noindent \it Proof. \rm}

\newtheorem{Theorem}{Theorem}
\newtheorem{Lemma}{Lemma}[section]

\newtheorem{Corollary}{Corollary}[section]
\newtheorem{Remark}{Remark}[section]

\@addtoreset{equation}{section}

\begin{document}

\title{Quasi-periodic solutions  for p-Laplacian equations with jumping nonlinearity and unbounded potential
terms
\author{Xiao Ma \footnote{  Email: mx5657@sina.cn}
, \ \ Daxiong  \ Piao \footnote{ The corresponding
author.Supported by the NSF of Shangdong
Province(No.ZR2012AM018),\ \ E-mail: dxpiao@ouc.edu.cn}
\\
        School of Mathematical Sciences,\ Ocean University of China,\ \ P.R. China
\\
        Yiqian Wang \footnote{Supported by the NNSF of Chian(No.11271183),Email:yqwangnju@yahoo.com}
\\
      Department of Mathematics,\ Nanjing University
\\
        \ \ Nanjing 210093,\ \ P.R. China
        \\
        \\
}
\date{}}

\maketitle

\begin{abstract}In this paper, we are concerned with the boundedness of all the solutions for a kind of second
order differential equations with p-Laplacian term $(\phi_p(x'))'+a\phi_p(x^+)-b\phi_p(x^-)+f(x)=e(t)$,
where $x^+=\max (x,0)$, $x^- =\max(-x,0)$, $\phi_p(s)=|s|^{p-2}s$, $p\geq2$, $a$ and $b$ are positive constants $(a\not=b)$,
and satisfy $\frac{1}{a^{\frac{1}{p}}}+\frac{1}{b^{\frac{1}{p}}}=2\omega^{-1} $,where $\omega \in \RR^+ \backslash \QQ$,
the perturbation $f$ is unbounded, $e(t)\in {\cal C}^{6}$ is is a smooth $2\pi_p$-periodic function on $t$,
where $\pi_p=\frac{2\pi(p-1)^{\frac{1}{p}}}{p\sin\frac{\pi}{p}}$.
\end{abstract}

\vskip 1.0cm
\section{ Introduction }

   Due to the relevance with applied mechanics, for example  modelling  some kind of suspension bridge(see\cite{[LazerM]}),  the  following semilinear Duffing's
   equations is widely studied
\begin {equation}
\label {tag1.1}
x^{\prime\prime} +ax^+ -bx^- =f(x, t),
\end {equation}
where $x^+=\max (x, 0),\  x^- =\max(-x, 0),\ f(x, t) $ is a smooth
$2\pi$-periodic function on $t$, $a $ and $b$ are positive
constants $(a\not=b)$.

If $f(x,t)$ depends only on $t$, the equation (\ref{tag1.1}) becomes
\beq\label{tag1.01}
x^{\prime\prime} +ax^+ -bx^- =f(t),\quad f(t+2\pi)=f(t),
\eeq
which had been studied by Fucik \cite{[Fucik]} and Dancer \cite{[Dancer]} in their
investigations of boundary value problems associated to  equations with
``jumping nonlinearities''. For recent developments, we refer to
\cite{[Gall], Habets, [Lazer]} and references therein.

  In 1996, Ortega \cite{[O-London]} proved the Lagrangian stability for the equation
\begin{equation}
\label {tag1.2}
  x'' +ax^+ -bx^- =1+ \gamma h(t)
\end{equation}
if $|\gamma | $ is sufficiently small and $h\in {\cal
C}^4(\mathbb{S}^1)$.

  On the other hand, when $\frac {1}{\sqrt{a}} +\frac {1}{\sqrt{b}}\in \QQ$,
Alonso and Ortega \cite{[Or]} proved that there is a $2\pi$-periodic function
$f(t)$ such that all the solutions of  Eq.~(\ref{tag1.01}) with large
initial conditions are unbounded. Moreover for such a $f(t)$, Eq.~(\ref{tag1.01})
has periodic solutions.

In 1999, Liu \cite{Liujmaa}  removed the smallness assumption on $|\gamma |$ in
  Eq.~(\ref{tag1.2}) when $\frac{1}{\sqrt {a}} +\frac {1}{\sqrt{b}} \in \QQ $ and obtained the same result.

Liu \cite {Liu1} studied the boundedness and Aubry-Mather sets for
the semilinear equation
$$
x''+\lambda^2x+\phi(x)=e(t),
$$
where $\phi(x)=o(x)$ as $|x|\rightarrow +\infty$.The methods and
formulations developed in \cite {Liu1} are benefit to treat the
more general equation

\beq x''+ax^+-bx^-+\phi(x)= e(t). \eeq

For the equation,Wang \cite{[Wangx]} and Wang \cite{[Wangy]}
considered the Lagrangian stability when the perturbation
$\phi(x)$ is bounded. While Yuan \cite{[Yuanx]} considered the
existence of quasiperiodic solutions and Lagrangian stability when
$\phi(x)$ is unbounded.

Fabry and Mawhin  \cite{[FM]} investigated the equation
\beq\label{semi} x''+ax^+-bx^-=f(x)+g(x)+e(t) \eeq under some
appropriate conditions, they get the boundedness of all solutions.

Yang \cite{[Yangx2]} considered  more complicated nonlinear equation with p-Laplacian operator
\beq\label{semi}
((\phi_p(x'))'+(p-1)[a\phi_p(x^+)-b\phi_p(x^-)]+f(x)+g(x)=e(t).
\eeq

Using Moser's small twist theorem, he proved that all the solutions are bounded,when $\frac{1}{a^{\frac{1}{p}}}+\frac{1}{b^{\frac{1}{p}}}=\frac{2m}{n}$, $m,n\in \NN $, the perturbation $f(x)$ and the oscillating term $g$ are bounded. For the case when $\frac{1}{a^{\frac{1}{p}}}+\frac{1}{b^{\frac{1}{p}}}=2\omega^{-1} $, where $\omega \in \RR^+ \backslash \QQ$, the perturbation $f(x)$ is bounded, Yang \cite{[Yangx1]} studied the following equation
\beq
(\phi_p(x'))'+a\phi_p(x^+)-b\phi_p(x^-)+f(x)=e(t).
\eeq
and came to the conclusion that every solution of the equation is bounded.

In 2004, Liu \cite{[Liujde04]} studied equation
\beq
(\phi_p(x'))'+a\phi_p(x^+)-b\phi_p(x^-)=f(x,t),f(x,t+2\pi)=f(x,t)
\eeq
where $p>1$, for the cases when $\frac{\pi_p}{a^{\frac{1}{p}}}+\frac{\pi_p}{b^{\frac{1}{p}}}=\frac{2\pi}{n}$
and $f\in{\cal C}^{(7,6)}(\RR\times \RR/2\pi\ZZ)$ and satisfies that

   {\rm \ (i)}\  the following limits exists {\rm uniformly\ in\ } $t$
  $$
  \  \lim_{x\rightarrow \infty}f(x,t)=f_{\pm}(t)\quad
$$

{\rm (ii)}\ the following limits exists {\rm uniformly\ in\ } $t$
$$
\ \lim_{x\rightarrow \infty}x^m \frac{\pa^{m+n}}{\pa x^m\pa t^n}f(x,t)=f_{\pm,m,n}(t)
$$
for $(n,m)=(0,6),\ (7,0)$ and $(7, 6)$. Moveover, $f_{\pm,m,n}(t)\equiv 0$ for $m=6,\ n=0, 7$. He comes to the conclusion that all solutions are bounded and the existence of quasi-periodic solutions.\\

In 2012,Jiao,Piao and Wang \cite{[JPW]} considered the
bounededness of equations \beq
x''+\omega^2x+\phi(x)=G_x(x,t)+p(t), \eeq and \beq
x''+ax^+-bx^-=G_x(x,t)+p(t). \eeq

Inspired by the above references, especially by Liu\cite {Liu1}
and Yuan\cite{[Yuanx]},  we are going to study the boundedness of
all solutions for the following equation \beq \label{tag(1.10)}
(\phi_p(x'))'+a\phi_p(x^+)-b\phi_p(x^-)+f(x)=e(t) \eeq where
$f(x)$ is unbounded, which generalized the equation in
\cite{[Yuanx]} for $p\geq2$.
 Our main results are as follows:\\
\begin{Theorem}  Assume\\
(H1) $f\in{\cal \CC}^{6}(\RR\backslash \{0\})\bigcap{\cal
\CC}^{0}(\RR)$, and there are constants $c$ and
$0<\gamma<\frac{1}{p-1}$ such that
\[|x^kf^{(k)}(x)|\leq c|x|^\gamma\]
for all $x\in \RR\setminus \{0\}$, where $ 0\leq k \leq6$;\\
(H2) There are two constants $\beta_1,\beta_2>0$, and
$\beta_1,\beta_2$ satisfy that $p\beta_1>q\beta_2>0$, where
$\frac{1}{p}+\frac{1}{q}=1$ such that for all $x\in \RR\setminus
\{0\}$, we have
\[xf(x)\geq\beta_1|x|^{\gamma+1},   x^2f'(x)\leq\beta_2|x|^{\gamma+1};\]\\
(H3) $e(t)\in{\cal \CC}^{6}(\mathbb{S}_p)$, where $\mathbb{S}_p\triangleq\RR/2\pi_p \ZZ$.\\

Then there exists $\delta_0>0$ such that for a given
$0<\rho<\frac{1}{2}$ and any $0<\delta<\delta_0$ and every
irrational number $\omega=\omega(\delta)$ satisfying
\beq\label{tag(1.12)}
1+\rho\leq\frac{\omega-\omega_0}{\delta}\leq2-\rho, \eeq where
$\omega_0=-\pi_p(\frac{1}{a^{\frac{1}{p}}}+\frac{1}{b^{\frac{1}{p}}}),
$ and \beq\label{tag(1.13)}
|m\omega-2\pi_pn|\geq\rho\delta|m|^{-\frac{3}{2}} \eeq for all
integers $n$ and $m\neq0$, the time $2\pi_p$ mapping $P$:
$(x,x')_{t=0}\rightarrow(x,x')_{t=2\pi_p}$ of the flow of
Eq.(\ref{tag(1.10)}) possesses an invariant curve $\Lambda_\delta$
with rotation number $\omega=\omega(\delta)$ and the curve
surrounds the origin $(x,x')=(0,0)$ and goes arbitrarily far from
the origin as $\delta\rightarrow0$. Moreover, the curve is an
intersection of the invariant torus in the $(x,x',t(mod
2\pi_p))$-space with the ${t=0}$-plane and any motion starting
from the torus is quasiperiodic , of basic frequencies $2\pi_p$
and $\omega$.
\end{Theorem}
\begin{Corollary}
The equation (\ref{tag(1.10)}) possesses Lagrange stability, i.e. if $x(t)$ is any solution of equation (\ref{tag(1.10)}), then it exists for all $t\in \RR$ and $sup_{t\in \RR}(|x(t)|+|x'(t)|)<\infty$.\\
\end{Corollary}
\begin{Corollary}
Most motions with large amplitude are quasiperiodic, i.e. most initial conditions (in the sense of Lebesgue measure) with large $|x(0)|+|x'(0)|$ give rise to quasiperidic motions: $x(t)=f(2\pi_p t,\omega t)$, where $f$ is a function on a 2-torus.
\end{Corollary}
\begin{Remark}
For the equation in (\ref{tag(1.10)}),We believe that for the
following two cases:(i) $\omega\in\QQ$;(ii) $1 <p<2$ similar
results still hold true. We will study those problems in our
future work.
\end{Remark}

The main idea is as follows: By means of transformation theory the original system outside of a
large disc $D=\{(x, x')\in \RR^{2}:  x^{2}+x'^{2}\leq r^{2} \}$ in
$(x, x')$-plane is transformed into a perturbation of an integrable
Hamiltonian system. The Poincar\'{e} map of the transformed system
is closed to a so-called twist map in $\RR^{2}\backslash D$. Then
Moser's twist theorem guarantees the existence of arbitrarily large
invariant curves diffeomorphic to circles and surrounding the origin
in the $(x,x')$-plane. Every such curve is the base of a
time-periodic and flow-invariant cylinder in the extended phase
space $(x,x',t)\in \RR^{2}\times \RR,$ which confines the solutions in
the interior and which leads to a bound of these solutions.\\

The remain part of this paper is organized as follows. In section 2,
we introduce action-angle variables and exchange the role of time and angle variables.  And we construct
canonical transformations such that the new Hamiltonian system is
closed to an integrable one. In section 3, we give some estimates. In section 4, we will prove the
Theorem 1 and the Corollaries by Moser's twist theorem.

Throughout this paper, $F(x)=\int_0^x f(s)ds,F(0)=0$, $c$ and $C$  are some positive constants without concerning their quantity.

\section{Action-Angle Variables }

In this section,  we will introduce action-angle variables
$(r,\theta)$ via symplectic transformations.

We introduce a new variables $y$ as $y=-\varphi_p(\omega^{-1} x)$,
let $q$ be the conjugate exponent of $p$ : $p^{-1}+q^{-1}=1$. Then
(\ref{tag(1.10)}) is changed into the form \beq
x'=-\omega\varphi_q(y),
y'=\omega[a_1\varphi_p(x^{+})-b_1\varphi_p(x^{-})]+\omega^{1-p}[f(x)-e(t)]
\eeq where $a=\omega^{p}a_1$, $b=\omega^{p}b_1$ and $a_1$, $b_1$
satisfy \beq a_1^{-\frac{1}{p}}+b_1^{-\frac{1}{p}}=2, \eeq which
is a planar non-autonomous Hamiltonian system \beq x'=-\frac{\pa
H}{\pa y }(x,y,t),  y'=\frac{\pa H}{\pa y }(x,y,t) \eeq where
$$H(x,y,t)=\frac{\omega}{q}|y|^q+\frac{\omega}{p}(a_1|x^+|^p+b_1|x^-|^p)+\omega^{1-p}[F(x)-e(t)x].
$$ Let $C(t)=\sin_p t$ be the solution of the following initial
value problem \beq (\varphi_p(C'(t)))'+\varphi_p(C(t))=0, \quad
C(0)=0,C'(0)=1. \eeq
  Then it follows from [16] that $C(t)=\sin_p(t)$ is a $2\pi_p$-period $C^2$ odd function with $\sin_p(\pi_p-t)=\sin_p(t)$, for $t\in[0,\frac{\pi_p}{2}]$ and
$\sin_p(2\pi_p-t)=-\sin_p(t)$, for $t\in[\pi_p,2\pi_p]$. Moreover for $t\in(0,\frac{\pi_p}{2})$, $C(t)>0$, $C'(t)>0$, and $C: [0,\frac{\pi_p}{2}]\rightarrow[0,(p-1)^{\frac{1}{p}}]$ can be implicitly given by
$$\int_0^{\sin_p t}\frac{ds}{(1-\frac{s^p}{p-1})^{\frac{1}{p}}}=t. $$
\begin{Lemma}
For $p\geq2$ and for any $(x_0,y_0)\in \RR^2$, $t_0\in \RR$, the solution $$z(t)=(x(t, t_0, x_0, y_0),y(t, t_0, x_0, y_0))$$ of (2.1) satisfying the initial condition $z(t_0)=(x_0,y_0)$ is unique and exists on the whole $t$-axis.
\end{Lemma}
 The proof of uniqueness can be obtained similarly as the proof of Proposition 2 in [17], the global existence result can be proved similarly as Lemma 3.1 in [10].\
Consider an auxiliary equation
\[(\phi_p(x'))'+a_1\phi_p(x^+)-b_1\phi_p(x^-)=0\]
Let $v(t)$ be the solution with initial condition: $(v(0),v'(0))=((p-1)^{\frac{1}{p}},0)$. Setting $\phi_p(v')=u $, then $(v,u)$ is a solution of the following planar system:
\[x'=\phi_q(y),\quad y'=-a_1\phi_p(x^+)+b_1\phi_p(x^-)\]
where $q=p/(p-1)>1$. It is not difficult to prove that:\\
(i) $q^{-1}|u|^q+p^{-1}(a_1|v^+|^p+b_1|v^-|^p)\equiv\frac{a_1}{q}$;\\
(ii) $v(t)$ and $u(t)$ are $2\pi_p$-periodic functions.\\
(iii)$v(t)$ can be given by
\beq
v(t)=\left\{
\begin{array}{ll}
\sin_p(a_1^{\frac{1}{p}}t+\frac{\pi_p}{2}),&\mbox{$0\leq t\leq\frac{\pi_p}{2a_1^{\frac{1}{p}}}$,}\\
-(\frac{a_1}{b_1})^{\frac{1}{p}}\sin_p{b_1}^{\frac{1}{p}}(t-\frac{\pi_p}{2a_1^{\frac{1}{p}}}),& \mbox{$
\frac{\pi_p}{2a_1^{\frac{1}{p}}}<t\leq\pi_p$.}
\end{array}
\right.
\eeq
\beq
v(2\pi_p-t)=v(t), t\in[\pi_p,2\pi_p].
\eeq
\begin{Lemma}
Let $I_p=\int_0^{\frac{\pi_p}{2}}\sin_p tdt.$ Then $$I_p=\frac{(p-1)^{\frac{2}{p}}}{p}B(\frac{2}{p},1-\frac{1}{p}),$$
where $B(r,s)=\int_0^1t^{r-1}(1-t)^{s-1}dt$ for $r>0,s>0.$
\end{Lemma}
 From the expression of $v(t)$ in (2.5), we obtain
\beq
\int_0^{\frac{\pi_p}{2a_1^{\frac{1}{p}}}}v(t)dt=\frac{I_p}{a_1^{\frac{1}{p}}},
\eeq
\beq
\int_{\frac{\pi_p}{2a_1^{\frac{1}{p}}}}^{\pi_p}v(t)dt=-\frac{a_1^{\frac{1}{p}}I_p}{b_1^{\frac{2}{p}}}.
\eeq
This method has been used in \cite{[Yangx1]}.

We introduce the action and angle variables via the solution
$(v(t),u(t))$ as follows. \beq\label{tag(2.9)}
x=d^{\frac{1}{p}}r^{\frac{1}{p}}v(\theta),y=d^{\frac{1}{q}}r^{\frac{1}{q}}u(\theta)
\eeq where $d=qa_1^{-1}$. This transformation is called a
generalized symplectic transformation as its Jacobian is 1. Under
this transformation, the system (2.1) is changed to \beq
\theta'=\frac{\partial h}{\partial
r}(r,\theta,t),r'=-\frac{\partial h}{\partial \theta}(r,\theta,t)
\eeq with the Hamiltonian function \beq \label{tag(2.10)}
h(r,\theta,t)=\omega
r+\omega^{1-p}F(d^{\frac{1}{p}}r^{\frac{1}{p}}v(\theta))-\omega^{1-p}d^{\frac{1}{p}}r^{\frac{1}{p}}v(\theta)e(t)
\eeq Observing the facts $F(\cdot)\in{\cal C}^{7}(\RR\backslash
\{0\})\bigcap{\cal C}^{1}(\RR)$ by (H1) and $e(t)\in{\cal
C}^{6}(\mathbb{S}_p)$ by (A3) and $v(\theta)\in{\cal
C}^{1}(\mathbb{S}_p)$, we have that $h(r,\theta,t)\in{\cal
C}^{1,1,6}(\RR\times \mathbb{S}_p \times \mathbb{S}_p)$. Let
$\Xi=\{\theta\in \mathbb{S}_p:v(\theta)=0\}$, clearly, the
Lebesgue measure of the set $\Xi$ vanishes. From the above, we
have that $h(r,\theta,t)$ is of class ${\cal C}^{6}$ in $r$ when
$\theta\in \mathbb{S}_p\setminus \Xi$ and $t\in \mathbb{S}_p$ are
regarded as parameters. Note that
\[rd\theta-hdt=-(hdt-rd\theta)\]
This means that if we can solve $r=r(h,t,\theta)$ from (\ref{tag(2.10)}) as a function of $h$, $t$ and $\theta$, then $r=r(h,t,\theta)$ is the Hamiltonian function of the following system:
\beq\label{tag(2.11)}
\frac{dh}{d\theta}=-\frac{\partial r}{\partial t}(h,t,\theta),\frac{dt}{d\theta}=\frac{\partial r}{\partial h}(h,t,\theta)
\eeq
i.e. (\ref{tag(2.11)}) is a Hamiltonian system with $r=r(h,t,\theta)$ as Hamiltonian function. Now the new action, angle and time variables are $h$, $t$, $\theta$. This method have been used in Levi \cite{LE1}.

\section{Some Estimates}
In this section, we will give some estimates which will be
used in the proof of Theorem. Throughout this section,  we
assume the hypotheses of Theorem 1 hold.
\begin{Lemma}
For $r$ large enough and $t\in \mathbb{S}_p$, it holds that:
\beq
 |\frac{\partial ^k }{\partial r^k }F(d^{\frac{1}{p}}r^{\frac{1}{p}}v(\theta)) |\leq c\cdot r^{-k+\frac{\gamma+1}{p}}, \quad0\leq k \leq 6
\eeq
\beq
|\frac{\partial ^k }{\partial r^k }f(d^{\frac{1}{p}}r^{\frac{1}{p}}v(\theta)) |\leq c\cdot r^{-k+\frac{\gamma+1}{p}}, \quad0\leq k \leq 6
\eeq
where $\theta\in \mathbb{S}_p$ if $k=1$ or $\theta\in \mathbb{S}_p\setminus\Xi$ if $k\geq2$.\\
\Proof Using the assumption (H1) and letting
$x=d^{\frac{1}{p}}r^{\frac{1}{p}}v(\theta)$, we have
\[|F(d^{\frac{1}{p}}r^{\frac{1}{p}}v(\theta))|=|\int_0^x f(s)ds|\leq\int_0^{|x|}|f(s)|ds\leq c|x|^{\gamma+1}\leq cr^{\frac{\gamma+1}{p}}\]
This means that (3.1) holds for $k=0$. Observe that
\[\frac{\partial ^k }{\partial r^k }F(d^{\frac{1}{p}}r^{\frac{1}{p}}v(\theta))=\sum_{s=1}^ k c_sF^{(s)}(d^{\frac{1}{p}}r^{\frac{1}{p}}v(\theta))\partial_r^{j_1}(d^{\frac{1}{p}}r^{\frac{1}{p}}v(\theta))\cdots\partial_r^{j_s}(d^{\frac{1}{p}}r^{\frac{1}{p}}v(\theta))\]
where $1\leq j_1,\cdots,j_s \leq k$, $j_1+\cdots+j_s = k$. Hence the terms are bounded by
\[c|F^{(s)}(d^{\frac{1}{p}}r^{\frac{1}{p}}v(\theta))||\partial_r^{j_1}(d^{\frac{1}{p}}r^{\frac{1}{p}}v(\theta))|\cdots|\partial_r^{j_s}(d^{\frac{1}{p}}r^{\frac{1}{p}}v(\theta))|\]
\[\leq c|F^{(s)}(d^{\frac{1}{p}}r^{\frac{1}{p}}v(\theta))||d^{\frac{1}{p}}r^{-j_1+\frac{1}{p}}v(\theta)|\cdots|d^{\frac{1}{p}}r^{-j_s+\frac{1}{p}}v(\theta)|\]
\[=c|f^{(s-1)}(d^{\frac{1}{p}}r^{\frac{1}{p}}v(\theta))(d^{\frac{1}{p}}r^{\frac{1}{p}}v(\theta))^{s-1}||d^{\frac{1}{p}}r^{\frac{1}{p}}v(\theta)|r^{-j_1-\cdots-j_s}\]
\[\leq c|r^{\frac{1}{p}}v(\theta)|^\gamma|r^{\frac{1}{p}}v(\theta)|r^{-k}\leq cr^{-k+\frac{\gamma+1}{p}},\]
where $\theta\in \mathbb{S}_p$ if $k=1$ or $\theta\in \mathbb{S}_p\setminus\Xi$ if $k\geq2$. This ends the proof of (3.1). The proof of (3.2) is similar to the one above.
\end{Lemma}

\vskip 3mm

\begin{Lemma}\label{lemmaJ}  Let
\beq\label{tag(3.3)}
h_1(r,\theta,t)=\omega^{1-p}F(d^{\frac{1}{p}}r^{\frac{1}{p}}v(\theta))-\omega^{1-p}d^{\frac{1}{p}}r^{\frac{1}{p}}v(\theta)e(t).
\eeq
Then, for $r$ large enough and $t\in \mathbb{S}_p$, we have
\begin {equation}
|\partial_r^k\partial_t^lh_1(r,\theta,t)|\leq cr^{-k+\frac{\gamma+1}{p}},\quad0\leq k \leq 6,
\end {equation}
where $\theta\in \mathbb{S}_p$ if $k=1$ or $\theta\in \mathbb{S}_p\setminus\Xi$ if $k\geq2$. \\
\Proof
The proof is finished by using Lemma 3.1 and (A3).\\

It follows from (\ref{tag(2.10)}) and (\ref{tag(3.3)}) that
\beq
h(r,\theta,t)=\omega r + h_1(r,\theta,t).
\eeq
Using Lemma 3.2 and noting $0<\gamma<\frac{1}{p-1}\leq1$,we have that, for $r$ large enough and $t\in \mathbb{S}_p$, the following inequalities hold:
\beq\label{tag(3.6)}
0<c_1r\leq h(r,\theta,t)\leq c_2r,\quad \theta\in \mathbb{S}_p,
\eeq
\beq
\partial_r h(r,\theta,t)\geq\frac{1}{2}\omega>0,\quad \theta\in \mathbb{S}_p,
\eeq
\beq
|\partial_r^k\partial_t^lh(r,\theta,t)|\leq cr^{-k+1}\leq cr^{-k}h(r,\theta,t),\quad 0\leq k+l\leq6,\quad\theta\in \mathbb{S}_p\setminus\Xi.
\eeq
Using the implicit theorem and (3.7) we can, indeed, solve (\ref{tag(2.10)}) for $r=r(h,t,\theta)$ as a function of $h$, $t$ and $\theta$.
\end{Lemma}

In the following, we will give Lemma 3.3, which is Lemma A1.1 in
M.Levi \cite{LE1}, for convenience, here we just give the Lemma
and omit the proof.
\begin{Lemma}
If a real function $f$ of two real variables $x$, $t$ ($t$ viewed as a parameter) satisfies for some $c>0$ and $n\in \NN$:
\[|\partial_x^k\partial_t^if(x,t)|\leq cx^{-k}f(x,t)\]
for all $x>0$ large enough and for all $k$,$i$:$k+i\leq \NN$ and if, moreover,
\[\partial_x f(x,t)\geq cx^{-1}f(x,t)>0\]
for all $x>0$ large enough, then the inverse function $g(y,t)$ of $f$ in $x$ satisfies
\[|\partial_y^k\partial_t^ig(y,t)|\leq cy^{-k}g(y,t),\]
for all $k+i\leq \NN$ and for all $y$ large enough.
\end{Lemma}

\begin{Lemma}
For $h$ large enough and $t\in \mathbb{S}_p$, it holds that:
\beq
 |\partial_h^k \partial_t^lr(h,t,\theta)|\leq ch^{-k}r(h,t,\theta)\leq c h^{-k+1}, \quad0\leq k \leq 6,\quad\theta\in \mathbb{S}_p\setminus\Xi.
\eeq
\Proof
Regard $\theta$ as a parameter. The proof is finished by using Lemma 3.3 together with (3.6) to (3.8).
\end{Lemma}

Let
\beq
g(r,\theta,t):=r^{-\frac{1}{p}}h_1(r,\theta,t)=r^{-\frac{1}{p}}\omega^{1-p}[F(d^{\frac{1}{p}}r^{\frac{1}{p}}v(\theta))-d^{\frac{1}{p}}r^{\frac{1}{p}}v(\theta)e(t)].
\eeq
By Lemma 3.2 we have that for $r\gg1$ and $\theta,t\in \mathbb{S}_p$,
\beq
|\partial_r^k\partial_t^lg(r,\theta,t)|\leq cr^{-k+\frac{\gamma}{2}},\quad0\leq k+l\leq6.
\eeq

Now (3.5) can be rewritten in the form of
\beq
h=\omega r+r^{\frac{1}{p}}g(r,\theta,t), \quad where \quad r=r(h,t,\theta).
\eeq
Hence, by Taylor's formula,
\[g(r,\theta,t)=g(\omega^{-1}h-\omega^{-1}r^{\frac{1}{p}}g(r,\theta,t),\theta,t)\]
\[=g(\omega^{-1}h,\theta,t)-\int_0^1g_r^{'}(\omega^{-1}h-s\omega^{-1}r^{\frac{1}{p}}g(r,\theta,t),\theta,t)\omega^{-1}r^{\frac{1}{p}}g(r,\theta,t)ds\]
\[=g(\omega^{-1}h,\theta,t)+R_0(h,t,\theta)\]
\beq
=(\omega^{-1}h)^{-\frac{1}{p}}\omega^{1-p}[F(d^{\frac{1}{p}}(\omega^{-1}h)^{\frac{1}{p}}v(\omega^{-1}\theta))-d^{\frac{1}{p}}(\omega^{-1}h)^{\frac{1}{p}}v(\omega^{-1}\theta)e(t)]+R_0(h,t,\theta)
\eeq
where
\beq
R_0(h,t,\theta)=-\int_0^1g_r^{'}(\omega^{-1}h-s\omega^{-1}r^{\frac{1}{p}}g(r,\theta,t),\theta,t)\omega^{-1}r^{\frac{1}{p}}g(r,\theta,t)ds
\eeq
By (3.12) we have
\[r=\omega^{-1}h-\omega^{-1}r^{\frac{1}{p}}g(r,\theta,t)\]
\[=\omega^{-1}h-\omega^{-1}(\omega^{-1}h-\omega^{-1}r^{\frac{1}{p}}g(r,\theta,t))^{\frac{1}{p}}g\]
\[=\omega^{-1}h-\omega^{-1}(\omega^{-1}h)^{\frac{1}{p}}(1-h^{-1}r^{\frac{1}{p}}g(r,\theta,t))^{\frac{1}{p}}g\]
\beq
=\omega^{-1}h-\omega^{-1}(\omega^{-1}h)^{\frac{1}{p}}g+\frac{1}{p}\omega^{-1}(\omega^{-1}h)^{\frac{1}{p}}g\int_0^1(1-sh^{-1}r^{\frac{1}{p}}g(r,\theta,t))^{\frac{1}{p}-1}h^{-1}r^{\frac{1}{p}}gds
\eeq
where we have expressed $(1-h^{-1}r^{\frac{1}{p}}g(r,\theta,t))^{\frac{1}{p}}$ by Taylor's formula. Putting (3.13) into the second term in the last line of (3.15) we have
\beq
r(h,t,\theta)=\omega^{-1}h-\omega^{-p}F(d^{\frac{1}{p}}(\omega^{-1}h)^{\frac{1}{p}}v(\omega^{-1}\theta))+R_1+R_2+R_3,
\eeq
where
\beq
 R_1=\omega^{-(2+\frac{1}{p})}h^{\frac{1}{p}}\int_0^1g_r^{'}(\omega^{-1}h-s\omega^{-1}r^{\frac{1}{p}}g(r,\theta,t),\theta,t)r^{\frac{1}{p}}g(r,\theta,t)ds
\eeq
\beq
R_2=\frac{1}{p}\omega^{-(1+\frac{1}{p})}h^{\frac{1}{p}-1}\int_0^1(1-sh^{-1}r^{\frac{1}{p}}g(r,\theta,t))^{\frac{1}{p}-1}r^{\frac{1}{p}}g^2ds
\eeq
\beq
R_3=d^{\frac{1}{p}}\omega^{-(p+\frac{1}{p})}h^{\frac{1}{p}}v(\theta)e(t)
\eeq
where $r=r(h,t,\theta)$, $g=g(r(h,t,\theta),\theta,t)$.
\begin{Lemma}
For $h\gg1$ , $t\in \mathbb{S}_p$, $\theta\in \mathbb{S}_p\setminus\Xi$, we have
\beq
|\partial_h^k\partial_t^lR_1(h,t,\theta)|\leq ch^{-k+\gamma},\quad 0\leq k+l\leq5,
\eeq
\beq
|\partial_h^k\partial_t^lR_2(h,t,\theta)|\leq ch^{-k+\gamma},\quad 0\leq k+l\leq5,
\eeq
\beq
|\partial_h^k\partial_t^lR_3(h,t,\theta)|\leq ch^{-k+\frac{1}{p}},\quad 0\leq k+l\leq5.
\eeq
\Proof
By (3.9) and (3.6), it is easy ti verify that for $s\in[0,1]$, $h\gg1$,  $t\in \mathbb{S}_p$, and $\theta\in \mathbb{S}_p\setminus\Xi$,
\beq
|\partial_h^k\partial_t^lr^{\frac{1}{p}}(h,t,\theta)|\leq ch^{-k+\frac{1}{p}},\quad 0\leq k+l\leq6
\eeq
Noting(3.11) and (3.23) we have that for any $s\in[0,1]$,
\beq
|\partial_h^k\partial_t^l(s\omega^{-1}r^{\frac{1}{p}}g)|\leq ch^{-k+\frac{1+\gamma}{p}},\quad 0\leq k+l\leq6
\eeq
Applying $\partial_h^k\partial_t^l$ to $g_r^{'}(\omega^{-1}h-s\omega^{-1}r^{\frac{1}{p}}g(r,\theta,t),\theta,t)$ with $r=r(h,t,\theta)$ and using (3.11), (3.24) together with $0<\gamma<\frac{1}{p-1}\leq1$, we have
\beq
|\partial_h^k\partial_t^lg_r^{'}(\omega^{-1}h-s\omega^{-1}r^{\frac{1}{p}}g(r,\theta,t),\theta,t)|\leq ch^{-k-1+\frac{\gamma}{p}},\quad 0\leq k+l\leq5.
\eeq
By applying $\partial_h^k\partial_t^l$ to $gh^{\frac{1}{p}}r^{\frac{1}{p}}$ and noting (3.23) and (3.11), we obtain
\beq
|\partial_h^k\partial_t^l(gh^{\frac{1}{p}}r^{\frac{1}{p}})|\leq ch^{-k+\frac{2}{p}+\frac{\gamma}{p}},\quad 0\leq k+l\leq6.
\eeq
By combining (3.17), (3.25), (3.26) and $p\geq2$, the proof of (3.20) is completed. The proof of (3.21) is similarly completed by using (3.9) and (3.11). The proof of (3.22) is obvious.
\end{Lemma}

\begin{Lemma}
Let \beq
\bar{F}(h)=\int_0^{2\pi_p}F(d^{\frac{1}{p}}(\omega^{-1}h)^{\frac{1}{p}}v(\theta))d\theta.
\eeq For $h\gg1$, we have the estimates: \beq
|\bar{F}^{(k)}(h)|\leq ch^{-k+\frac{\gamma+1}{p}},\quad 0\leq
k\leq6, \eeq

\beq \bar{F}^{'}(h)\geq ch^{-1+\frac{\gamma+1}{p}},
\eeq \beq
\bar{F}^{''}(h)\leq -ch^{-2+\frac{\gamma+1}{p}}.
\eeq
\Proof
It follows from (3.1) that (3.28) holds true. For simplicity we write $x=d^{\frac{1}{p}}(\omega^{-1}h)^{\frac{1}{p}}v(\theta)$.\\
Using (3.27) and the fact that the set $\Xi\bigcap[0,2\pi_p]$ is finite, we have
\[\bar{F}^{'}(h)=\int_{[0,2\pi_p]\setminus\Xi}F'(d^{\frac{1}{p}}(\omega^{-1}h)^{\frac{1}{p}}v(\theta))\frac{1}{p}d^{\frac{1}{p}}(\omega^{-1}h)^{\frac{1}{p}-1}\omega^{-1}v(\theta)d\theta\]
\beq =\frac{1}{ph}\int_{[0,2\pi_p]\setminus\Xi} x f(x)d\theta.
\eeq In view of assumption (H2), we have \beq
\bar{F}^{'}(h)=\frac{1}{ph}\int_{[0,2\pi_p]\setminus\Xi} x
f(x)d\theta\geq\frac{\beta_1}{ph}\int_{[0,2\pi_p]\setminus\Xi}
|x|^{\gamma+1}d\theta=ch^{-1+\frac{\gamma+1}{p}}. \eeq
This ends the proof of (3.29).\\

Differentiating (3.31) with respect to $h$ we have
\[\bar{F}^{''}(h)=\frac{1}{p^2h^2}\int_{[0,2\pi_p]\setminus\Xi} x^2 f(x)d\theta-\frac{1}{qh}\bar{F}^{'}(h)\]
\[\leq\frac{\beta_2}{p^2h^2}\int_{[0,2\pi_p]\setminus\Xi}|x|^{\gamma+1}d\theta-\frac{1}{qh}\bar{F}^{'}(h)\]
\[\leq\frac{\beta_1}{ph}\frac{\beta_2}{\beta_1}\bar{F}^{'}(h)-\frac{1}{qh}\bar{F}^{'}(h)\]
\[\leq(\frac{1}{p}\frac{\beta_2}{\beta_1}-\frac{1}{q})\frac{1}{h}\bar{F}^{'}(h)\]
\[\leq \frac{1}{p}(\frac{\beta_2}{\beta_1}-\frac{p}{q})ch^{-2+\frac{\gamma+1}{p}}.\]
The proof of (3.30) is now completed by assumption (H2):
$\beta_1>\beta_2>0$, $p\beta_1>q\beta_2>0$ and $p\geq 2$.
\end{Lemma}

Let
\beq
S(h,\theta)=\omega^{-p}\int_0^\theta(F(d^{\frac{1}{p}}(\omega^{-1}h)^{\frac{1}{p}}v(\vartheta))-\bar{F}(h))d\vartheta.
\eeq
In view of Lemma 3.1 and 3.6, we have
\beq
|\partial_h^{k+1}S(h,\theta)|\leq ch^{-k-1+\frac{\gamma+1}{p}}.
\eeq
Define a map $\Psi_1$ : $(h,t)\rightarrow(h,\tau)$ by formula:
\[\Psi_1\quad :\quad h=h,\quad t=\tau+\frac{\partial S}{\partial h},\]
where the time variable $\theta$ is regarded as a parameter. Clearly the map is of period $2\pi_p$ in $\theta$. Since $dh\wedge dt=dh\wedge d\tau$, the map is symplectic. The Hamiltonian function $r=r(h,t,\theta)$ defined in (3.16) is transformed into a new Hamiltonian function
\[\hat{r}=\omega^{-1}h-\omega^{-p}F(d^{\frac{1}{p}}(\omega^{-1}h)^{\frac{1}{p}}v(\theta))+R_1+R_2+R_3+\frac{\partial S}{\partial \theta},\]
i.e.
\beq
\hat{r}=r(h,\tau+\partial_hS(h,\theta),\theta)=\omega^{-1}h-\omega^{-p}\bar{F}(h)+R(h,\tau,\theta),
\eeq
where
\beq
R(h,\tau,\theta):=(R_1+R_2+R_3)(h,\tau+\partial_hS(h,\theta),\theta).
\eeq
Applying $\partial_h^k\partial_\tau^l$ to (3.36) and using Lemma 3.4 together with (3.34), we have that for $h\gg1$, $\tau\in \mathbb{S}_p$, $\theta\in \mathbb{S}_p\setminus\Xi$,
\beq
|\partial_h^k\partial_\tau^lR(h,\tau,\theta)|\leq ch^{-k+\max{\{\gamma,\frac{1}{p}}\}},\quad0\leq k+l\leq5.
\eeq
This method have been used in [13].

\section{The proof of theorem}
For given $0<\delta<1$, define a map $\Psi_2$ : $(h,\tau)\rightarrow(\lambda,\tau)$ by $\Psi_2$ : $\delta\lambda=\omega^{-p}\bar{F}'(h)$, $\tau=\tau$, $1\leq\lambda\leq4$. This trick was used in \cite{LE1}. Observe that small $\delta$ corresponds to large $h$ since $\bar{F}'(h)\rightarrow0$ as $h\rightarrow+\infty$ by (3.28) with $k=1$. It follows from (3.30) that $\Psi_2$ is a diffeomorphism. The equations corresponding to $\hat{r}(h,\tau,-\theta)$ are transformed into
\beq
\frac{d\lambda}{d\theta}=l_1(\lambda,\tau,\theta,\delta),\quad\frac{d\tau}{d\theta}=-\omega^{-1}+\delta\lambda+l_2(\lambda,\tau,\theta,\delta),
\eeq
where
\beq
l_1(\lambda,\tau,\theta,\delta)=\omega^{-p}\delta^{-1}\bar{F}''(h)\partial_\tau R(h,\tau,-\theta),
\eeq
\beq
l_2(\lambda,\tau,\theta,\delta)=-\partial_h R(h,\tau,-\theta),
\eeq
with $h=h(\delta\lambda)$ being defined by $\delta\lambda=\omega^{-p}\bar{F}'(h)$.
\begin{Lemma}
For $h(\delta\lambda)$ we have
\beq
c_1\delta^{\frac{p}{\gamma+1-p}}\leq h(\delta\lambda)\leq c_2\delta^{\frac{p}{\gamma+1-p}},
\eeq
\beq
|\partial_\lambda^kh(\delta\lambda)|\leq c h(\delta\lambda),\quad0\leq k\leq4.
\eeq
\Proof
The abbreviation $h(\delta\lambda)=h$ is used in what follows. By (3.28) and (3.29) we have
\beq
c_3h^{-1+\frac{\gamma+1}{p}}\leq\bar{F}'(h)\leq c_4h^{-1+\frac{\gamma+1}{p}}.
\eeq
Putting $h(\delta\lambda)$ into (4.6) and observing $\bar{F}'(h)=\omega^{p}\delta\lambda$, we have
\[c_3h^{-1+\frac{\gamma+1}{p}}\leq\omega^{p}\delta\lambda\leq c_4h^{-1+\frac{\gamma+1}{p}},\]
by direct computing we have (4.4). In the following, we give the proof of (4.5), which is similar to the reference \cite{LE1}. For convenience, here we just give a brief proof.\\

Differentiating $\bar{F}'(h)=\omega^{p}\delta\lambda$ with respect to $\lambda$ we have
\beq
\partial_\lambda h\cdot\bar{F}''(h)=\omega^{p}\delta,
\eeq \beq (\partial_\lambda
h)^2\cdot\bar{F}'''(h)+\partial_\lambda^2 h\cdot\bar{F}''(h)=0
\eeq \beq (\partial_\lambda
h)^3\cdot\bar{F}^{(4)}(h)+3\partial_\lambda
h\cdot\partial_\lambda^2h\cdot\bar{F}'''(h)+\partial_\lambda^3
h\cdot\bar{F}''(h)=0,
 \eeq
 From (4.7) and (3.30) we have
\[|\partial_\lambda h|=|\frac{\omega^{p}\delta}{\bar{F}''(h)}|=|\frac{\omega^{p}\delta h}{h\bar{F}''(h)}|=c|\frac{\delta h}{\delta\lambda}|\leq ch,\]
By (4.8) and (3.30) we have
\[|\partial_\lambda^2h|=|-\frac{\bar{F}'''(h)}{\bar{F}''(h)}(\partial_\lambda h)^2|\leq ch,\]
similarly we can get
\[|\partial_\lambda^kh|\leq ch.\]
Now we finish the proof of the Lemma.
\end{Lemma}
\begin{Lemma}
For $0<\delta\ll1$, $\tau\in \mathbb{S}_p$, $\theta\in \mathbb{S}_p\setminus\Xi$, $\lambda\in[1,4]$, we have
\beq
|\partial_\lambda^k\partial_\tau^ll_1(\lambda,\tau,\theta,\delta)|,\quad|\partial_\lambda^k\partial_\tau^ll_2(\lambda,\tau,\theta,\delta)|\leq c\delta^\sigma,\quad0\leq k+l\leq4,
\eeq
where $\sigma=(-1+\max{\{\gamma,\frac{1}{p}\}})\frac{p}{\gamma+1-p}>1$.
\Proof
Since $l_2(\lambda,\tau,\theta,\delta)=-\partial_h R(h,\tau,-\theta)$, by the estimate (3.37) we have
\[|\partial_\tau^ll_2|=|\partial_h\partial_\tau^lR|\leq ch^{-1+\max{\{\gamma,\frac{1}{p}\}}}\leq c\delta^{(-1+\max{\{\gamma,\frac{1}{p}\}})\frac{p}{\gamma+1-p}}=c\delta^\sigma.\]
This shows that (4.10) holds for $l_2$ and $k=0$ and $0\leq
l\leq4$. When $k>0$ and $\theta\in \mathbb{S}_p\setminus\Xi$, by
application of $\partial_\lambda^k\partial_\tau^l$ to both sides
of $l_2(\lambda,\tau,\theta,\delta)=-\partial_h R(h,\tau,-\theta)$
we have
\[\frac{\partial^{k+l}}{\partial_\lambda^k\partial_\tau^l}l_2=\sum_{s=1}^k\frac{\partial^{s+1+l}}{\partial_\lambda^{s+1}\partial_\tau^l}R(h,\tau,\theta)\cdot(\frac{\partial}{\partial\lambda})^{j_1}h(\frac{\partial}{\partial\lambda})^{j_2}h\cdots (\frac{\partial}{\partial\lambda})^{j_s}h,\]
where $1\leq j_1,j_2,\cdots,j_s\leq k$, $j_1+j_2+\cdots+j_s=k$. Consequently, using (3.37) and Lemma 4.1 we have
\[\left|\frac{\partial^{k+l}}{\partial_\lambda^k\partial_\tau^l}l_2\right|\leq\sum_{s=1}^{k}ch^{-s-1+\max{\{\gamma,\frac{1}{p}\}}}\cdot h^s\leq ch^{-1+\max{\{\gamma,\frac{1}{p}\}}}\leq c\delta^\sigma.\]
This completes the proof of (4.10) for $l_2$.\\

Using the same trick above and using Lemma 4.1 together with (3.28), (3.37) and (4.2), the remaining proof can be finished.
\end{Lemma}
\begin{Lemma}
The solutions $(\lambda(\theta),\tau(\theta))$ of  Eq.(4.1) with the initial conditions $\lambda(0)=\lambda_0\in[2,3]$, $\tau(0)=\tau_0$ do exist for $0\leq\theta\leq2\pi_p$ if the $\delta$ is sufficiently small. The Poncar\'{e} mapping of  Eq.(4.1) is of the form
\beq\label{latestsystem}{\cal{P}}^{2\pi_p}:\quad \left\{\begin{array}{ll}
\lambda(2\pi_p)=\lambda_0+\delta L_1(\lambda_0,\tau_0,\delta),\\
\tau(2\pi_p)=\tau_0-2\pi_p\omega^{-1}+\delta (\lambda_0+L_2(\lambda_0,\tau_0,\delta)),
\end{array}\right.
\eeq
where $L_1$, $L_2$ satisfy
\beq
|\partial_{\lambda_0}^k\partial_{\tau_0}^lL_1(\lambda_0,\tau_0,\delta)|,\quad|\partial_{\lambda_0}^k\partial_{\tau_0}^lL_2(\lambda_0,\tau_0,\delta)|\leq c\delta^{\sigma-1},\quad0\leq k+l\leq4.
\eeq
\Proof
By integrating Eq.(4.1) from $\theta=0$ to $\theta=2\pi_p$, and using Lemma 4.2 and the contraction principle, the proof can be easily finished. The argument  is similar to the one of Lemma 4 in \cite{[DZ]}. We omit the details.
\end{Lemma}
\vskip 0.3cm
\noindent {\it Proof of Theorem 1 and the Corollaries.}
Since  $l_1$, $l_2$ satisfy (4.10), it is easy to verify that the solutions $(\lambda(\theta),\tau(\theta))$ of  Eq.(4.1) with the initial conditions $\lambda(0)=\lambda_0\in[2,3]$, $\tau(0)=\tau_0$ do exist for $0\leq\theta\leq2\pi_p$ if the $\delta$ is sufficiently small,that is the conditions of Lemma 4.3 come true. Hence the poincr\'{e} map in the form of (4.11) does exist, and the map has the intersection property in the domain $[2,3]\times \mathbb{S}_p$, i.e. if $\Gamma$ is an embedded circle in $[2,3]\times \mathbb{S}_p$ homotopic to a circle $\lambda=const$, then $P(\Gamma)\bigcap{\Gamma}\neq\varnothing$. This is a well-known fact. See \cite{[DZ]}, for instance. \\

Until now we have verified that the mapping $P$ satisfy all the conditions of Moser's small twist theorem \cite{[Mos]} in the domain $[2,3]\times \mathbb{S}_p$, if  $\delta$ is small enough. We come to the conclusion that for any $0<\delta<\delta_0$ with some constants $\delta_0$ small enough, the mapping $P$ has an invariant curve $\Upsilon_\delta$ in the annulus $[2,3]\times \mathbb{S}_p$ with rotation number $\omega=\omega(\delta)$ satisfying (\ref{tag(1.12)}) and (\ref{tag(1.13)}). Retracting the sequence of the transformations back to the original Eq.(\ref{tag(1.10)}), we conclude that the time $2\pi_p$ mapping $P$: $(x,x')_{t=0}\rightarrow(x,x')_{t=2\pi_p}$ of the flow of Eq. (\ref{tag(1.10)}) possesses an invariant curve $\Gamma_\delta$ which surrounds the origin in the $(x,y)$-plane. Going back to (\ref{tag(2.10)}) and using the fact that small $\delta$ corresponds to large $h$ and inequality (\ref{tag(3.6)}) we know that $r=r(h,t,\theta)\rightarrow+\infty$ as $\delta\rightarrow 0$. Returning to (\ref{tag(2.9)}) and using the formula (i) in the second section, we have $x^2+y^2\rightarrow\infty$ as $\delta\rightarrow0$. Thus the invariant curve $\Gamma_\delta$ goes arbitrarily far from the origin in the $(x,y)$-plane when $\delta\rightarrow0$. Since Eq. (\ref{tag(1.10)})is of period $2\pi_p$ in time $t$, these curves $\Lambda_\delta$ are the intersections of the invariant tori in the $(x,x',t(mod 2\pi_p))$-space with the ${t=0}$-plane and any motion stating from the torus is quasiperiodic of basic frequencies $(2\pi_p,\omega)$. This ends the proof of the Theorem 1. As there exist invariant curves of the poncar\'{e} mapping of the system (\ref{tag(1.10)}), which surround the origin $(x,y)=(0,0)$ and are arbitrarily far from the origin. According to the Moser's small twist theorem , we can know that the system (\ref{tag(1.10)}) possesses Lagrange stability. The statement of the first Corollary has been proved. And according to \cite{LE1} and \cite{[DZ]} , we can get the second Corollary.

\end{document}